\documentclass[reqno]{amsart}
\usepackage{graphicx}
\usepackage{amssymb, amsmath, amsthm}
\usepackage{amsfonts}
\usepackage{amssymb}
\usepackage{float}
\usepackage{mathrsfs}
\usepackage{enumitem} 
\newcommand{\dis}{\displaystyle}

\newcommand{\divv}{\text{\rm div}}

\newcommand{\R}{\mathbb R}

\numberwithin{equation}{section}
\numberwithin{theorem}{section}

\numberwithin{figure}{section}

\begin{document}

\title[Nonexistence of Self-similar Singularities]{Nonexistence of Self-similar Singularities in the Ideal Viscoelastic Flow}
 
\author{Anthony Suen} 

\address{Department of Mathematics\\Indiana University Bloomington, IN 47405}

\email{cksuen@indiana.edu}

\dedicatory{To my family and my wife Candy}

\date{October 25, 2011} 


\keywords{viscoelastic flow; self-similar singularities} 

\subjclass[2000]{35A20} 

\begin{abstract}
We prove the nonexistence of finite time self-similar singularities in an ideal viscoelastic flow in $\R^3$. We exclude the occurrence of Leray-type self-similar singularities under suitable integrability conditions on velocity and deformation tensor. We also prove the nonexistence of asymptotically self-similar singularities in our system. The present work extends the results obtained by Chae in the case of magnetohydrodynamics (MHD).
\end{abstract}

\maketitle

\section{Introduction}
We prove the nonexistence of finite time self-similar singularity in the Cauchy problem of the ideal viscoelastic flow in $\R^3$:
\begin{align}
\label{1.1} \left\{ \begin{array}{l}
u_t + (u\cdot\nabla)u=-\nabla p+\nabla\cdot FF^T,\\
F^k_t+(u\cdot\nabla)F^k=(F^k\cdot\nabla)u,\\
\divv(u)=0,\\
u(x,0)=u_0,F(x,0)=F_0(x).
\end{array}\right.
\end{align}
Here $x\in\R^3$ is the spatial variable and $t\ge0$ is the time, $u=u(x,t)=(u^1,u^2 ,u^3)$ is the velocity of the flow, $P=P(x,t)$ is a prescribed scalar pressure, $F=F(x,t)\in\R^{3\times3}$ is the deformation tensor with $F^k$ being the k-th column of $F$. The initial data $(u_0,F_0)$ is assumed to satisfy the following condition:
\begin{align}\label{1.2}
\divv(u_0)=\divv(F^k_0)=0.
\end{align}
By taking divergence on the second equation in \eqref{1.1} and using $\divv(u)=0$,
\begin{align*}
\divv(F^k)_t+(u\cdot\nabla)(\divv(F^k))=0,
\end{align*} 
so $F$ is ``divergence free'' for later time. In view of the above, we can rewrite \eqref{1.1}-\eqref{1.2} into the following equivalent form:
\begin{align}
\label{1.4} \left\{ \begin{array}{l}
u_t + (u\cdot\nabla)u=-\nabla p+\sum(F^k\cdot\nabla)F^k,\\
F^k_t+(u\cdot\nabla)F^k=(F^k\cdot\nabla)u,\\
\divv(u)=\divv(F^k)=0,\\
u(x,0)=u_0,F(x,0)=F_0(x).
\end{array}\right.
\end{align}

We refer to \cite{bcah} or \cite{dafermos} for a more detailed derivation and physical discussion of the above system.

The study of incompressible fluids can be dated back to 18th century when Euler \cite{euler} considered the well-known Euler's equations:
\begin{align}\label{1.3}
 \left\{ \begin{array}{l}
u_t + (u\cdot\nabla)u=-\nabla p,\\
\divv(u)=0.
\end{array}\right.
\end{align}
The local-in-time existence of \eqref{1.3} can be found in Kato \cite{kato}, yet it is still interesting to know whether or not finite time singularities do occur in the above system. There seems to be no satisfactory result until Beale et.al.\cite{bkm} discovered that, if a solution $u$ for \eqref{1.3} possesses a singularity at a finite time $T>0$, then it is necessary to have
\begin{align*}
\int_0^T||\nabla\times u(\cdot,t)||_{L^\infty}dt=+\infty.
\end{align*}
The above blow-up criterion is later generalized to many other related systems, for example ideal magnetohydrodynamics (MHD) (Caflisch-Klapper-Steel \cite{cks}) and ideal viscoelastic flow (Hu and Hynd \cite{hu}). With a clever argument suggested by Chae in \cite{chae1}-\cite{chae2}, he applied the results in Caflisch et.al.\cite{cks} to exclude the possibility of a finite time apparition of self-similar singularities in both viscous and ideal MHD model. The goal of our present work is to extend Chae's results to ideal viscoelastic flow based on the blow-up criterion shown by Hu and Hynd \cite{hu}.

\medskip

We now give a precise formulation of our result. To begin with, we observe that if $(u,F,p)$ is a solution of \eqref{1.4}, then for any $\lambda>0$, $(u^{(\lambda)},F^{(\lambda)},p^{(\lambda)})$ is also a solution with the initial data $(u_0^{(\lambda)},F_0^{(\lambda)})$, where
\begin{align*}
u^{(\lambda)}(x,t)&=\lambda u(\lambda x,\lambda^2t),\\
F^{(\lambda)}(x,t)&=\lambda F(\lambda x,\lambda^2t),\\
p^{(\lambda)}(x,t)&=\lambda^2 p(\lambda x,\lambda^2t),\\
(u_0^{(\lambda)}(x),F_0^{(\lambda)}(x))&=(\lambda u_0(\lambda x),\lambda F_0(\lambda x)).
\end{align*}
Therefore, if there exists a self-similar blowing up solution $(u(x,t),F(x,t))$ of \eqref{1.4}, then it has to be of the form
\begin{align}\label{1.6}
u(x,t)=\frac{1}{\sqrt{T_{*}-t}}U\left(\frac{x}{\sqrt{T_{*}-t}}\right),
\end{align}
\begin{align}\label{1.7}
F^k(x,t)=\frac{1}{\sqrt{T_{*}-t}}E^k\left(\frac{x}{\sqrt{T_{*}-t}}\right),
\end{align}
\begin{align}\label{1.8}
p(x,t)=\frac{1}{T_{*}-t}P\left(\frac{x}{\sqrt{T_{*}-t}}\right),
\end{align}
when $t$ is being close to the possible blow-up time $T_*$. By substituting \eqref{1.6}-\eqref{1.8} into \eqref{1.4}, we find that $(U,E,P)=(U(y),E(y),P(y))$ satisfies of the following stationary system:
\begin{align}
\label{1.9} \left\{ \begin{array}{l}
\frac{1}{2}U+\frac{1}{2}(y\cdot\nabla)U+(U\cdot\nabla)U=-\nabla P+\sum\limits_{k}(E^k\cdot\nabla)E^k,\\
\frac{1}{2}E^k+\frac{1}{2}(y\cdot\nabla)E^k+(U\cdot\nabla)E^k=(E^k\cdot\nabla)U,\\
\divv(U)=\divv(E^k)=0.
\end{array}\right.
\end{align}

Conversely, if $(U,E,P)$ is a smooth solution of \eqref{1.9}, then $(u,F,p)$ as defined by \eqref{1.6}-\eqref{1.8} is a classical solution of \eqref{1.4} which blows up at $t=T_*$. Leray \cite{leray} was the first one to study self-similar singularities of the form similar to \eqref{1.6}-\eqref{1.8} in Navier-Stokes equations, while Chae \cite{chae1}-\cite{chae2} considered those in the case of (MHD). We apply those concepts to viscoelastic flow which will be given in later sections.

\medskip

The following Theorem~1.1-1.3 are the main results of this paper:

\medskip

\noindent{\bf Theorem~1.1} \em Suppose there exists $T_{*}>0$ such that we have a representation of a solution $(u,F)$ to \eqref{1.4} by \eqref{1.6}-\eqref{1.8} for all $t\in(0,T_*)$ with $(U,E)$ satisfying the following conditions:
\begin{enumerate}[label={\upshape{(1.\arabic*)}}, ref={\upshape{2.}\arabic*}, topsep=0.3cm, itemsep=0.15cm] 
\setcounter{enumi}{8}
\item $(U,E)\in(C_0^1(\R^3))^2$,
\item there exists $q>0$ such that $(\Omega,E)\in (L^{r}(\R^3))^2$ for all $r\in(0,q)$, where $\Omega=\nabla\times U$.
\end{enumerate}
Then we have $U=E=0$.
\rm
\medskip

Theorem~1.1 can be proved in a way exactly the same as given in Chae \cite{chae2} pp.1014--1017, so we omit the proof here. By changing the decay conditions on $(U,E)$, we can derive similar results about the nonexistence of self-similar singularities. The following Theorem~1.2 is reminiscent of Theorem~1.2 in \cite{chae1}:

\medskip

\noindent{\bf Theorem~1.2} \em Suppose there exists $T_{*}>0$ such that we have a representation of a solution $(u,F)$ to \eqref{1.4} by \eqref{1.6}-\eqref{1.8} for all $t\in(0,T_*)$ with $(U,E)$ satisfying the following conditions:
\begin{enumerate}[label={\upshape{(1.\arabic*)}}, ref={\upshape{2.}\arabic*}, topsep=0.3cm, itemsep=0.15cm] 
\setcounter{enumi}{10}
\item $(U,E)\in(H^m(\R^3))^2$ for $m>\frac{3}{2}+1$,
\item $||\nabla U||_{L^\infty}+||\nabla E||_{L^\infty}<\varepsilon$, \\where $\varepsilon>0$ is a constant as chosen in {\rm Theorem~2.1}.
\end{enumerate}
Then we have $U=E=0$.
\rm
\medskip

Using Theorem~1.1-1.2, we can obtain the following result about asymptotically self-similar singularity, which is reminiscent of Theorem~1.3 in \cite{chae1}. We refer to Chae \cite{chae2} (especially section~2) for more detailed descriptions of asymptotically self-similar singularities.
\medskip

\noindent{\bf Theorem~1.3} \em Given $T>0$ and $m>\frac{3}{2}+1$, let $(u,F)\in((C([0,T);H^m(\R^3)))^2$ be a classical solutions to \eqref{1.4}. Suppose there exists functions $U$, $E$ satisfying {\rm(1.8)-(1.9)} as in {\rm Theorem~1.1} such that
\setcounter{equation}{12}
\begin{align}\label{1.12}
&\sup_{0<t<T}\frac{1}{T-t}\left|\left| u(\cdot,t)-\frac{1}{\sqrt{T-t}}U\left(\frac{\cdot}{\sqrt{T-t}}\right) \right|\right|_{L^1}\\
&\qquad+\sup_{0<t<T}\frac{1}{T-t}\left|\left| F(\cdot,t)-\frac{1}{\sqrt{T-t}}E\left(\frac{\cdot}{\sqrt{T-t}}\right) \right|\right|_{L^1}<\infty,\notag
\end{align}
and
\begin{align}\label{1.13}
&\lim\limits_{t\nearrow T}(T-t)\left|\left|\nabla u(\cdot,t)-\frac{1}{\sqrt{T-t}}\nabla U\left(\frac{\cdot}{\sqrt{T-t}}\right) \right|\right|_{L^\infty}\\
&\qquad+\lim\limits_{t\nearrow T}(T-t)\left|\left|\nabla F(\cdot,t)-\frac{1}{\sqrt{T-t}}\nabla E\left(\frac{\cdot}{\sqrt{T-t}}\right) \right|\right|_{L^\infty}=0,\notag
\end{align}
then we have $U=E=0$. 

Moreover, there exists $\delta>0$ such that $(u,F)$ can be extended to a solution of \eqref{1.4} in $[0,T+\delta]\times\R^3$ with $(u,F)\in C([0,T+\delta];H^m(\R^3))$.
\rm

\bigskip

This paper is organized as follows. In section~2 we first prove Theorem~2.1 about a continuation principle for local solution of \eqref{1.4} with the help of an auxiliary lemma. We then begin the proofs of Theorem~1.2-1.3 in section~3 which basically follow the arguments as given in Chae \cite{chae1}-\cite{chae2}.

\medskip

\section{Continuation Principle for Local Solution}

In this section we show the continuation principle for local solution of \eqref{1.4} with the help of an energy estimate as given in Lemma~2.2. The result is reminiscent of Lemma~2.1 in Chae \cite{chae1} except that we apply the blow-up criterion for singularities in ideal viscoelastic flow proved by Hu and Hynd \cite{hu}. We begin with the following Theorem~2.1:
\medskip

\noindent{\bf Theorem~2.1} \em Given $T>0$ and $m>\frac{3}{2}+1$, let $(u,F)\in((C([0,T);H^m(\R^3)))^2$ be a classical solutions to \eqref{1.4}. There exists $\varepsilon>0$ and $\delta>0$ such that if 
\begin{align}\label{2.1}
\sup_{0\le t<T}(T-t)\left[||\nabla u(\cdot,t)||_{L^\infty}+||\nabla F(\cdot,t)||_{L^\infty}\right]<\varepsilon,
\end{align}
then $(u,F)$ can be extended to a solution of \eqref{1.4} in $[0,T+\delta]\times\R^3$ satisfying $(u,F)\in C([0,T+\delta];H^m(\R^3))$.
\rm
\medskip

Theorem~2.1 can be proved by an argument developed by Chae \cite{chae1}. It requires an energy estimate which is given by the following lemma:
\medskip

\noindent{\bf Lemma~2.2} \em Assume that the hypotheses and notations as in {\rm Theorem~2.1} are in force. Then for $0\le t<T$, we have the following estimates
\begin{align}\label{2.2}
&\frac{d}{dt}\left(||u||^2_{H^m}+\sum_k||F^k||^2_{H^m}\right)\\
&\qquad\qquad\le M\left(||\nabla u||_{L^\infty}+\sum_k||\nabla F^k||^2_{L^\infty}\right)\left(||u||^2_{H^m}+\sum_k||F^k||^2_{H^m}\right),\notag
\end{align}
where $M>0$ is a generic constant which depends only on $m$.
\begin{proof}
The proof is exactly as in Hu and Hynd \cite{hu}, pp.4--6.
\end{proof}
\rm
\medskip

\begin{proof}[\bf proof of Theorem~2.1] We follow the proof of Lemma~2.1 as in Chae \cite{chae1}. We claim
\begin{align}\label{2.3}
\int_0^T(||\nabla u||^2_{L^\infty}+\sum_k||\nabla F^k||^2_{L^\infty})ds<\infty,
\end{align}
and if \eqref{2.3} holds, then we have
\begin{align*}
\int_0^T(||\nabla\times u||_{L^\infty}+\sum_k||\nabla\times F^k||_{L^\infty})ds<\infty,
\end{align*}
and so by the blow-up criterion derived in \cite{hu}, there exists $\delta>0$ such that $(u,F)$ can be continuously extended to a solution of \eqref{1.4} in $[0,T+\delta]\times\R^3$ with $(u,F)\in (C([0,T+\delta];H^m(\R^3)))^2$.

To prove \eqref{2.3}, we apply \eqref{2.2} as in Lemma~2.2 to obtain
\begin{align}\label{2.4}
&\frac{d}{dt}\left[(T-t)(||u||^2_{H^m}+\sum_k||F^k||^2_{H^m})\right]+||u||^2_{H^m}+\sum_k||F^k||^2_{H^m}\\
&\qquad\le M(T-t)(||\nabla u||_{L^\infty}+\sum_k||F^k||_{L^\infty})(||u||^2_{H^m}+\sum_k||F^k||^2_{H^m}).\notag
\end{align}
If we choose $\varepsilon<\frac{1}{2M}$, then \eqref{2.1} and \eqref{2.4} imply
\begin{align*}
\frac{d}{dt}\left[(T-t)(||u||^2_{H^m}+\sum_k||F^k||^2_{H^m})\right]+\frac{1}{2}\left(||u||^2_{H^m}+\sum_k||F^k||^2_{H^m}\right)\le0.
\end{align*}
Upon integrating the above,
\begin{align}\label{2.5}
\dis\int_0^T(||u||^2_{H^m}+&\sum_k||F^k||^2_{H^m})ds\le 2T(||u_0||^2_{H^m}+\sum_k||F_0^k||^2_{H^m})<\infty
\end{align}
By the Sobolev embedding theorem, there exists $l\ge1$ and $\beta\in(0,1)$ such that $H^m(\R^3)\hookrightarrow C^{l,\beta}(\R^3)$, so we obtain from \eqref{2.5} that
\begin{align*}
\int_0^T(||\nabla u||_{L^\infty}+\sum_k||\nabla F^k||_{L^\infty})ds\le\int_0^T(||u||_{H^m}+\sum_k||F^k||_{H^m})ds\\
\le\sqrt{T}\left[\int_{0}^T(||u||^2_{H^m}+\sum_k||F^k||^2_{H^m})ds\right]^\frac{1}{2}<\infty,
\end{align*}
and hence \eqref{2.3} follows.
\end{proof}
\medskip

\section{Nonexistence of Asymptotically Self-similar Singularities: \\Proof of Theorem~1.2-1.3}

In this section we prove Theorem~1.2-1.3 as stated in section~1. These arguments are reminiscent of Chae \cite{chae1} and hence we omit some of those technical details.
\medskip

\begin{proof}[\bf proof of Theorem~1.2]
By the definitions \eqref{1.6}-\eqref{1.7} for $t\in(0,T_{*})$,
\begin{align*}
||U(\cdot,t)||_{L^\infty}&=(T_{*}-t)||\nabla u(\cdot,t)||_{L^\infty},\\
||\nabla E(\cdot,t)||_{L^\infty}&=(T_{*}-t)||\nabla F(\cdot,t)||_{L^\infty}.
\end{align*}
So assumption {\rm (1.12)} implies $$\sup_{0\le t<T_{*}}(T_{*}-t)\left[||\nabla u(\cdot,t)||_{L^\infty}+||\nabla F(\cdot,t)||_{L^\infty}\right]<\varepsilon.$$
By Theorem~2.1, $(u,F)$ can be continuously extended to a solution of \eqref{1.4} in $[0,T_*+\delta]\times\R^3$, which is impossible unless $u=F=0$. Therefore $U=E=0$ as required.
\end{proof}
\medskip

\begin{proof}[\bf proof of Theorem~1.3]
Define $(\tilde u,\tilde F,\tilde p)$ by
\begin{align*}
\tilde u(y,s)&=(\sqrt{T-t})u(x,t),\\
\tilde F(y,s)&=(\sqrt{T-t})F(x,t),\\
\tilde p(y,s)&=(\sqrt{T-t})p(x,t),
\end{align*}
where $\dis y=\frac{x}{\sqrt{T-t}}$ and $\dis s=\frac{1}{2}\log\left(\frac{T}{T-t}\right)$. Assumptions \eqref{1.12}-\eqref{1.13} can then be rewritten as
\begin{align}
\sup_{0<s<\infty}||\tilde u(\cdot,s)-U(\cdot)||_{L^1}+\sup_{0<s<\infty}||\tilde F(\cdot,s)-E(\cdot)||_{L^1}<\infty
\end{align}
and
\begin{align}\label{3.2}
\lim\limits_{s\to\infty}||\nabla\tilde u(\cdot,s)-\nabla U(\cdot)||_{L^\infty}+||\nabla\tilde F(\cdot,s)-\nabla E(\cdot)||_{L^\infty}=0,
\end{align}
and we obtain $\lim\limits_{s\to\infty}||\tilde u(\cdot,t)-U(\cdot)||_{H^1(B_R)}=\lim\limits_{s\to\infty}||\tilde F(\cdot,t)-E(\cdot)||_{H^1(B_R)}=0$ for all $R>0$ where $B_R=\{x\in\R^3:|x|<R\}$. Hence $U,E\in C^1_0(\R^3)$. Apply the same argument as in Chae \cite{chae1} pp.451--453, we have that $(U,E)$ is a classical solution of \eqref{1.4} satisfying {\rm (1.10)}, and so Theorem~1.1 implies $U=E=0$. 

Next we put $U=E=0$ into \eqref{3.2}, then
\begin{align*}
\lim\limits_{s\to\infty}||\nabla\tilde u(\cdot,s)||_{L^\infty}+||\nabla\tilde F(\cdot,s)||_{L^\infty}=0,
\end{align*}
and so there exists $s'>0$ such that
$$\sup_{s'\le s<\infty}||\nabla\tilde u(\cdot,s)||_{L^\infty}+||\nabla\tilde F(\cdot,s)||_{L^\infty}<\varepsilon,$$
where $\varepsilon>0$ is as chosen in Theorem~2.1. Let $t'=T(1-e^{-2s'})$, then $(u,F)$ satisfies
$$\sup_{t'\le t<T}(T-t)(||\nabla u(\cdot,t)||_{L^\infty}+||\nabla F(\cdot,t)||_{L^\infty})<\varepsilon,$$
and by Theorem~2.1, there exists $\delta>0$ such that $(u,F)$ can be extended to a solution of \eqref{1.4} in $[0,T+\delta]\times\R^3$ with $(u,F)\in C([0,T+\delta];H^m(\R^3))$. This completes the proof of Theorem~1.3.
\end{proof}
\bigskip


\subsection*{Acknowledgment}
The author would like to thank Professor Dongho Chae for his comments and letting me to use his idea from his previous papers.

\end{document}